\documentclass[11pt,twoside]{amsart}

\usepackage{bbm,amsmath,amssymb,amsthm,amscd,latexsym,mathrsfs}
\usepackage[all]{xy}
\usepackage[pdftex]{hyperref}

\newtheorem{Prop}{Proposition}[section]
\newtheorem{Thm}[Prop]{Theorem}
\newtheorem*{Theorem}{Theorem}

\theoremstyle{definition}

\newcommand{\EE}{\mathcal{E}}

\newcommand{\MM}{\mathcal{M}}

\newcommand{\uu}{\mathbbmss{1}}

\newcommand{\CC}{\mathbb{C}}
\newcommand{\ZZ}{\mathbb{Z}}

\newcommand{\NN}{\mathbb{N}}
\newcommand{\QQ}{\mathbb{Q}}
\newcommand{\PP}{\mathbb{P}}
\newcommand{\XX}{\mathcal{X}}

\DeclareMathOperator{\Hilb}{Hilb}
\DeclareMathOperator{\ch}{ch}
\DeclareMathOperator{\cc}{c}

\begin{document}

\title{On the Andr\'e motive of certain irreducible symplectic varieties}
\author{Ulrich Schlickewei}

\address{Mathematisches Institut der Universit{\"a}t Bonn, Endenicher Allee 60,
53115 Bonn, Germany} 
\email{uli@math.uni-bonn.de}

\begin{abstract}
We show that if $Y$ is an algebraic deformation of the Hilbert square of a K3 surface, then
the Andr\'e motive of $Y$ is an object of the 
category generated be the motive of $Y$ truncated in degree 2.
\end{abstract}

\maketitle 
\let\thefootnote\relax\footnotetext{This paper 
is a part of my Ph.D.\ thesis prepared at the University of Bonn. I would like to thank my advisor
Daniel Huybrechts for his continuous encouragement. Moreover, I am grateful to Eyal Markman for
valuable comments on a previous version of this note.

This work was supported by the SFB/TR 45 `Periods,
Moduli Spaces and Arithmetic of Algebraic Varieties' of the DFG
(German Research Foundation) and by the Bonn International Graduate School in Mathematics (BIGS).}

An irreducible symplectic variety $Y$ is a smooth, projective variety over $\CC$ which is simply connected and
which admits a nowhere degenerate, holomorphic two-form $\sigma \in H^0(Y,\Omega^2_Y)$. A general principle
says that most of the geometry of $Y$ is encoded by the cohomology group $H^2(Y,\ZZ)$ together
with the Hodge decomposition and the Beauville--Bogomolov quadratic form. Beauville \cite{Be} found two series of examples of irreducible symplectic varieties. Apart from these only
two exceptional examples have been discovered by O'Grady \cite{OG1}, \cite{OG2}.

In this note, we study the Andr\'e motive of irreducible symplectic varieties which are deformation equivalent to
the Hilbert scheme of points on a K3 surface or equivalently which are deformations of a smooth, compact
moduli space of stable sheaves
on a K3 surface. In all even dimensions, there is one family of such deformations, 
these families build one of Beauville's series of examples. 
Let $Y$ be such a variety. Denote by $\mathfrak{h}^2(Y)$ 
the Andr\'e motive of $Y$ truncated in degree 2.
We use results of Markman
on the monodromy of moduli spaces of sheaves on K3 surfaces and Andr\'e's deformation principle to derive

\begin{Theorem}
  a)The motive of $Y$ is an object of $\langle \mathfrak{h}^2(Y) \rangle$, the category generated by $\mathfrak{h}^2(Y)$.

  b) The motive of $Y$ is an object of the category generated by motives of Abelian
  varieties.

  c) All Hodge classes on $Y$ are motivated, hence absolute in the sense of Deligne.
\end{Theorem}

Item a) can be seen as a motivic manifestation
of the above-mentioned principle. Items b) and c) are consequences of a) and of Andr\'e's results \cite{An2} on
degree 2 motives of
Hyperk\"ahler varieties.

It has been proved
by Arapura \cite{Ar} that the motive of a moduli space $Y$ parametrizing sheaves on a K3 surface $S$
is an object of the category generated by the motive
of $S$. Since the motive of $S$ is also an object of $\langle \mathfrak{h}^2(Y) \rangle$, our result
can be seen as a generalization of Arapura's result.

\vspace{2ex}

After quickly reviewing Andr\'e's motives in Section \ref{AndreMotives}, we collect in Section 
\ref{MarkmanResults} some of Markman's results \cite{Mar1}, \cite{Mar2} on the cohomology of moduli spaces of sheaves on K3 surfaces 
and on their monodromy groups. The proof of our result is given in Section \ref{Proof}.

\begin{section}{Andr\'e motives} \label{AndreMotives}
The idea of the category of motives is to provide the target of a universal cohomology functor
for smooth, projective varieties. Grothendieck dreamed of a category which should be Tannakian and semisimple. 
However, Jannsen \cite{Ja} proved that Grothendieck's category of homological motives
can only be semisimple if homological and numerical equivalence of algebraic
cycles coincide on all varieties. This is one of Grothendieck's standard conjectures
which are widely open.

In order to circumvent the standard conjectures and to obtain nonetheless a Tannakian and semisimple
category of motives, Andr\'e introduced a category in which he formally inverted the dual Lefschetz
operator. The basic ingredient in the theory is the notion of a \emph{motivated cohomology class}. 
Let $X$ be a smooth, projective variety over $\CC$. A cohomology class $\alpha \in H^*(X,\QQ)$ is 
motivated if there exist a smooth, projective variety $Y$ and 
algebraic cycles $Z_1, Z_2$ on $X \times Y$ such that 
\begin{equation*}
  \alpha = p_{1,*} \left( [Z_1] \cup \Lambda_{H \times G} [Z_2] \right).
\end{equation*}
Here, $p_1 : X \times Y \to X$ and $p_2: X \times Y \to Y$ are the projections, 
and $\Lambda_{H \times G}$ is the dual Lefschetz
operator with respect to some product polarization $p_1^*H + p_2^* G$.
Clearly, algebraic cohomology classes are motivated. Vice versa, Grothendieck's standard conjectures
would imply that motivated cohomology classes are algebraic.

Andr\'e \cite{An1} proves that the category of motives defined in terms of motivated correspondences is
Tannakian and semisimple over $\QQ$. The Betti realization which maps a motive to the underlying singular
cohomology group is a conservative fibre functor. (Recall that a functor $F: \mathcal{C} \to \mathcal{C'}$ 
is conservative if a morphism $f$ in $\mathcal{C}$ is an isomorphism if and only if so is $F(f)$.)

One of the big advantages of motivated cohomology classes is that examples are rather easy to produce. This is
mainly due to the following result which gives a positive answer to Grothendieck's invariant cycle
conjecture in the motivated world.

\begin{Thm}[Andr\'e, see \cite{An1}, 5.1] \label{AndreInvariantCycle}
Let $f: \XX \to S$ be a smooth projective morphism
where $S$ is a smooth, connected, algebraic variety. Let $s \in S$ be a closed point, $m,n \in \NN$ 
and let
\begin{equation*}
  \alpha \in H^* (\XX_s, \QQ)^{\otimes n} \otimes \big( H^*(\XX_s, \QQ)^{\vee} \big)^{\otimes m}
\end{equation*}
be a motivated class which is invariant under a subgroup of finite index of $\pi_1(S,s)$ (acting
on $H^*(S,\QQ)$ via the monodromy representation). Then any translate of $\alpha$ under parallel transport
to $H^* (\XX_t, \QQ)^{\otimes n} \otimes (H^*(\XX_t, \QQ)^{\vee})^{\otimes m}$ for $t \in S$ is motivated on $\XX_t$.
\end{Thm}

\end{section}

\begin{section}{Markman's results} \label{MarkmanResults}
Let $S$ be a projective K3 surface, polarized by an ample divisor $H$. 
The Todd genus of $S$ is $\mathrm{td}(S) = 1 + 2 [x]$ where $x$ is an 
arbitrary point of $S$. A square root is given by $\sqrt{\mathrm{td}(S)} = 1 + [x]$.

For a coherent sheaf $E$ on $S$ define the Mukai vector by 
\begin{equation*}
  v(E) = \ch(E) \sqrt{\mathrm{td}(S)} \in H^*(S,\ZZ).
\end{equation*}
We associate with $S$ a rational weight two Hodge structure
\begin{equation*}
 \widetilde{H}(S,\QQ) := H^*(S,\QQ), \; \; \widetilde{H}^{2,0}(S) = H^{2,0}(S), \; \; \widetilde{H}^{1,1}(S)
  = H^0(S) \oplus H^{1,1}(S) \oplus H^{4}(S).
\end{equation*}
There is a natural duality operator 
\begin{equation*}
  D_S : \widetilde{H}(S,\QQ) \to \widetilde{H}(S,\QQ)
\end{equation*}
acting as $(-1)^i \mathrm{id}$ on $H^{2i}(S,\QQ)$. Since the K\"unneth components of the
dia\-gonal in $S \times S$ are algebraic, $D_S$ is given by an algebraic class.

The Mukai pairing on $\widetilde{H}(S,\QQ)$ is given by
\begin{equation*}
 \langle \alpha, \beta \rangle  = - \int_S D_S (\alpha) \cup \beta.
\end{equation*}
This is a non-degenerate, symmetric bilinear form of signature $(4+ , 20-)$.

\vspace{2ex} 
Let now $v \in H^*(S,\ZZ)$ be a primitive and effective (cf.\ \cite[Def.\ 1.1]{Mar2}) vector. 
Then by results of Mukai, Huybrechts, O'Grady and Yoshioka, there exist a polarization $H$ on $S$ and a 
non-empty, smooth, projective variety $X := M_H(v)$ which parametrizes $H$-stable sheaves with Mukai vector $v$
on $S$. Moreover, $X$ is 
an irreducible symplectic variety of dimension $d =\langle v,v \rangle + 2$
which is deformation equivalent to $\Hilb^{\frac{d}{2}} (S)$. We assume that $d > 2$ and for simplicity
we assume that $X$ is a fine moduli space.

Let $\EE$ be a universal sheaf on $S \times X$. Then $\EE$
is uniquely determined up to the twist by the pull-back of a line bundle from $X$. Denote by 
$p: S \times X \to S$ and by $q: S \times X \to X$ the projections, let
$\pi_{H^2} : H^*(X,\QQ) \to H^2(X,\QQ)$ be the projection in degree 2. Define
\begin{equation*} 
   \varphi_1' : \widetilde{H}(S,\QQ) \to H^2(X,\QQ), \; \; \;
      \alpha \mapsto 
    \pi_{H^2} \big\{  q_* \big( \ch(\EE) \cup p^*\sqrt{\mathrm{td}(S)} \cup p^* D_S (\alpha) \big) \big\}.
\end{equation*}
According to a result of O'Grady \cite{OG}, the restriction of $\varphi'_1$ to $v^{\perp}$ is
an isomorphism of Hodge structures (even over $\ZZ$).

We normalize the correspondence $\ch(\EE) p^* \sqrt{ \mathrm{td}(S)}$
following \cite[Lemma 3.1]{Mar2}: let $\eta := \varphi'_1 (v) \langle v,v \rangle^{-1} \in H^2(X,\QQ)$.
Put 
\begin{equation} \label{Defu}
 u:= \ch(\EE) \cup p^* \sqrt{ \mathrm{td}(S)} \cup q^* \exp(\eta).
\end{equation}
Then $u$ is independent of the universal sheaf $\EE$ and we define
\begin{equation*} \varphi_1: 
 \widetilde{H}(S,\QQ) \to H^2(X,\QQ), \; \; \;
 \alpha \mapsto \pi_{H^2}  \big\{ q_* (u \cup p^* D_S (\alpha)) \big\}.
\end{equation*}
Note that 
\begin{equation*}
  \varphi_1 (\alpha) = \varphi_1' (\alpha) - \frac{\langle \alpha, v \rangle}{\langle v , v \rangle} 
   \varphi'_1 (v).
\end{equation*}
This implies that $\varphi_1(v) = 0$ and that
$\varphi_{1| v^{\perp}} = \varphi'_{1|v^{\perp}}$.

Next, we note that $\varphi_1$ is an algebraic correspondence. This is, because $u$ and $D_S$
are so and because the projection $H^*(X,\QQ) \to H^2(X,\QQ)$ is algebraic (cf.\ \cite{Ar}
where the conjecture $B$ is shown for $X$).

Since the standard conjecture $B$ holds for $S$ as well, there is an algebraic right
inverse $\psi: H^2(X,\QQ) \to \widetilde{H}(S,\QQ)$ (see \cite[Cor.\ 3.14]{Kl}).
Since the (Mukai-)orthogonal
projection $\widetilde{H}(S,\QQ) \to \QQ v$ is given on $S \times S$ by the class
$-(\langle v, v \rangle^{-1}  D_S(v)) \otimes v$, the orthogonal projection
$\widetilde{H}(S,\QQ) \to v^{\perp}$ is algebraic. Thus we may assume
that $\psi$
induces an isomorphism
\begin{equation} \label{DefPsi}
  \psi: H^2(X,\QQ) \stackrel{\sim}{\to} v^{\perp} \subset \widetilde{H}(S,\QQ)
\end{equation}
which is inverse to $\varphi_{1|v^{\perp}}$.

\vspace{2ex}
Let $G_v$ be the fix group of $v$ in $\mathrm{Aut}(\widetilde{H}(S,\ZZ), \langle \; , \; \rangle )$.
Markman defines two representations of $G_v$ on $H^*(X,\ZZ)$. We will now describe both
of them.

\vspace{2ex}
1.) Let $p_{ij}$ be the projection from $X \times S \times X$ to the $(i,j)$-th
factor. For $g \in G_v$ set
\begin{equation*}
 \gamma_g' := (p_{13})_* \Big( p_{12}^* D_{X \times S} \big( (\mathrm{id} \otimes g) ({}^t u) \big) \cup p_{23}^* u \Big)^{-1} \in H^*(X \times X, \QQ),
\end{equation*}
where $D_{X \times S}$ is the duality operator acting by $(-1)^i$ on $H^{2i}(X \times S,\QQ)$ and
the class $u$ was introduced in (\ref{Defu}).
Let $l: H^* (X,\QQ) \to H^*(X,\QQ)$ be the universal polynomial map which takes the Chern character 
$(r+ a_1 + a_2 + \ldots )$ of a coherent sheaf to its total Chern class $(1+ a_1 + (\frac{a_1^2}{2} - a_2) + \ldots)$.
Then by definition
\begin{equation*}
 \gamma_g := \text{degree} \; d \; \text{part of} \; l(\gamma'_g) = \cc_d(\gamma'_g).
\end{equation*}

\begin{Thm}[Markman, \cite{Mar2}, Thm.\ 3.10 and Cor.\ 3.14] \label{ThmMarkmanInvariance} 
i) For $g \in G_v$ the correspondence $\gamma_g$
acts as a (degree-preserving) automorphism on $H^*(X,\QQ)$.

\vspace{1ex}
ii) The map 
\begin{equation*}  \gamma:  G_v \to \mathrm{Aut}(H^*(X,\QQ)), \; \;
  g \mapsto  \gamma_g
\end{equation*}
is a faithful representation of $G_v$.

\vspace{1ex}
iii) The class $u \in \widetilde{H}(S,\QQ) \otimes H^*(X,\QQ)$ is invariant under the product representation 
of $G_v$, where $G_v$ acts on the first factor via the natural representation.
\end{Thm}

The theorem implies that the algebraic maps
\begin{equation} \label{varphi_i}  \varphi_i: 
  \widetilde{H}(S,\QQ) \to H^{2i}(X,\QQ), \; \;
  \alpha \mapsto \pi_{H^{2i}} \big\{ q_* (u \cup p^* D_S (\alpha)) \big\}
\end{equation}
are $G_v$-equivariant.

\vspace{2ex} 
2.) To define the second representation of $G_v$ recall that the space
$\widetilde{H}(S,\QQ)$ has four positive directions. Given any two positive four-spaces $F$ and $F'$,
orientations of these spaces can be compared using orthogonal projections. An isometry 
$g\in \mathrm{Aut}(\widetilde{H}(S,\QQ), \langle \; , \; \rangle)$ is called 
\emph{orientation preserving} if for an oriented, positive four-space $F$ the space  
$g(F)$ has the same orientation.
This induces the \emph{covariance} or \emph{orientation character} 
\begin{equation} \label{OrientationCharacter}
  \mathrm{cov}:G \to \ZZ / 2 \ZZ
\end{equation}
sending $g$ to $0$ or $1$ according
to whether it preserves orientations or not.

Then Markman defines the representation
\begin{equation*}  
 \gamma_{\mathrm{mon}}: G_v \to \mathrm{Aut} (H^*(X,\QQ)), \; \;
 g \mapsto (D_X)^{\mathrm{cov}(g)} \circ \gamma_g
\end{equation*}
where again $D_X$ is the duality operator of $X$, acting by $(-1)^i \mathrm{id}$
on $H^{2i}(X)$.

The subscript is justified by the following result of Markman. An ele\-ment $g
\in \mathrm{Aut}(H^*(X,\QQ))$ is called a \emph{monodromy operator} if there exist a family
$\XX \to B$ of Hyperk\"ahler manifolds 
with fibre $\XX_b = X$ for some $b\in B$ and a $\widetilde{g} \in \pi_1(B,b)$ such
that $g$ is the image of $\widetilde{g}$ under the monodromy representation 
$\pi_1(B,b) \to \mathrm{Aut}(H^*(X,\QQ))$. Let $\mathrm{Mon}(X)$ be the subgroup of
$\mathrm{Aut}(H^*(X,\QQ))$ generated by monodromy operators.

\begin{Thm}[Markman, \cite{Mar2}, Thm.\ 1.6] \label{ThmMarkmanMonodromy} 
 The image of the representation 
 $\gamma_{\mathrm{mon}} : G_v \to \mathrm{Aut}(H^*(X,\QQ))$ is a normal subgroup
 of finite index in $\mathrm{Mon}(X)$.
 
 In particular, since $\gamma$ and $\gamma_{\mathrm{mon}}$ coincide on the kernel $N'$
 of the orientation character, its image $N: = \gamma(N')$ in $\mathrm{Aut}(H^*(X,\QQ))$
 is a subgroup of finite index in $\mathrm{Mon}(X)$.
\end{Thm}

Following an idea of Beauville, Markman had proved in previous work that the class of
the diagonal in $X \times X$ can be expressed in terms of the Chern classes of the universal
sheaf $\EE$. This implies
\begin{Thm}[Markman, \cite{Mar1}, Cor.2] \label{ThmMarkmanKuenneth}
The K\"unneth factors of the Chern classes of the universal
sheaf $\EE$ generate $H^*(X,\QQ)$.
\end{Thm}

\noindent \emph{Summary of results used in the sequel.}
We have seen that there are 
homomorphisms $\varphi_i: \widetilde{H}(S,\QQ) \to
H^{2i}(X,\QQ)$ and $\psi: H^2(X,\QQ) \to v^{\perp}$ with the following properties:

i) The $\varphi_i$ and $\psi$ are induced by algebraic cycles on $S \times X$ resp.\ on $X \times S$.

ii) The homomorphism $\varphi_1$ induces an isomorphism 
 \begin{equation*}
   v^{\perp} \to H^2(X,\QQ)
 \end{equation*}
 whose inverse is $\psi$.
 
iii) There is a subgroup of finite index $N \subset \mathrm{Mon}(X)$ such that the compositions
 \begin{equation*}
   \eta_{i,0} := \varphi_i \circ \psi : H^2(X,\QQ) \to H^{2i}(X,\QQ)
  \end{equation*}
  are $N$-equivariant. This follows from Theorem \ref{ThmMarkmanInvariance}
  and \ref{ThmMarkmanMonodromy}.

 iv) For $i \ge 2$, the classes $\varphi_i (v) \in H^{2i}(X,\QQ)$ are $N$-invariant. Again, this
 is implied by Theorem \ref{ThmMarkmanInvariance}.

 v) The sum of $H^0(X,\QQ)$, of the image of $\oplus_{i \ge 1} \eta_{i,0}$
  and of the $\varphi_i (v)$ generate
  the cohomology ring $H^*(X,\QQ)$
  as a $\QQ$-algebra. This is a consequence of Theorem \ref{ThmMarkmanKuenneth}.
\end{section} 

\begin{section}{Proof of the Theorem} \label{Proof}
\noindent a) Let $X = M_H(v)$ be as in the previous section, let $Y$ be a fixed
algebraic deformation of $X$. By this we
mean that there exists a smooth, projective morphism of connected, smooth, complex algebraic varieties
$\XX \to B$ which admits $X$ and $Y$ as fibers.

We have to prove that $\mathfrak{h}(Y)$ is an object of $\langle \mathfrak{h}^2(Y) \rangle$.
Recall that by definition, this is the smallest full subcategory of
the category $\MM$ of Andr\'e motives 
which contains $\mathfrak{h}^2(Y)$ and the unit object $\uu = \mathfrak{h} (\mathrm{Spec}(\CC))$,
which is stable under $\otimes$ and under duals and which contains all subobjects resp.\ quotients
in $\MM$ of objects in $\langle \mathfrak{h}^2(Y) \rangle$.

\vspace{1ex}
 The idea is to identify $\widetilde{H}(S,\QQ)$ with $G(X) = H^2(X,\QQ) \oplus \QQ$ and to use
 Markman's results to define a surjection of a sum of products of $G(X)$ to $H^*(X,\QQ)$ which
 is monodromy invariant. By Andr\'e's deformation principle, this will induce
 a surjection of a motive $\mathfrak{m}(Y)$ to $\mathfrak{h}(Y)$ where $\mathfrak{m}(Y)$ is
 an object of $\langle \mathfrak{h}^2(Y) \rangle$. Let's make this precise now.

 For any fibre $V$ of $\XX \to B$, let $\mathfrak{g}_0(V) := \mathfrak{h}^0(V) \simeq \uu$
 and for $i=1, \ldots, d=\dim(X)$ define 
 \begin{equation*}
   \mathfrak{g}_i (V) := \big( \mathfrak{h}^2(V) \oplus \uu(-1) \big) (-i+1).
 \end{equation*}
 (For $V = X = M_H(v)$, the motive $\mathfrak{g}_i(X)$ plays the role of $\widetilde{\mathfrak{h}}(S)(-i+1) = 
 \big( \mathfrak{h}^0(S)(-1) \oplus
  \mathfrak{h}^2(S) \oplus \mathfrak{h}^4(S)(1)\big) (-i+1)$.)

 Next, we put
 \begin{equation*}
   \mathfrak{m}(V) := \bigoplus_{(i_1, \ldots, i_{d}) \in \{0, \ldots, d \}^{d} } \big( \mathfrak{g}_{i_1}(V) \otimes \ldots \otimes
    \mathfrak{g}_{i_d}(V) \big).
 \end{equation*}
 Note that $\mathfrak{m}(V)$ is an object of $\langle \mathfrak{h}^2(V) \rangle$ and that $\mathfrak{m}(V)$ can be seen as a submotive
 of the motive of a variety $Z(V)$ which is a disjoint union of products of $V$ and $\PP^1$.

 We fix an isomorphism $\eta_0: \uu \to \mathfrak{h}^0(X)$. For $i= 1, \ldots, d$ we will
 define below morphisms of motives
 \begin{equation*}
   \eta_i : \mathfrak{g}_i (X) \to \mathfrak{h}^{2i} (X)
 \end{equation*}
 with the following properties:

 \vspace{1ex}
 1) there exists a subgroup $N$ of finite index in $\mathrm{Mon}(X)$ such that $\eta_i$ is $N$-invariant.
 
 2) if we define the morphism 
  \begin{equation*}
    \eta : \mathfrak{m}(X) \to \mathfrak{h}(X)
  \end{equation*}
  as the composition of the morphism
  \begin{equation*}
   \bigoplus (\eta_{i_1} \otimes \ldots \otimes \eta_{i_d}): \mathfrak{m}(V) \to \bigoplus_{(i_1 \ldots, i_d) \in \{0, \ldots, d\}^d} \big( \mathfrak{h}^{2i_1}(X) \otimes
    \ldots \otimes \mathfrak{h}^{2 i_d}(X) \big) 
  \end{equation*}
  with the cup-product morphism
  \begin{equation*}
   \bigoplus \big( \mathfrak{h}^{2i_1}(X) \otimes \ldots \otimes \mathfrak{h}^{2i_d}(X) \big) \to \mathfrak{h}(X),
  \end{equation*}
  then $\eta$ is surjective.

  \vspace{1ex}
  Assume for one moment, that the $\eta_i$ are defined. Consider the family $\mathcal{Z} \to B$ which is constructed
  by letting vary the $Z(V)$ over $B$. The morphism $\eta$ corresponds to a motivated
  cohomology class on $Z(X) \times X$. Property 1) implies that this class is invariant under a subgroup of
  finite index of the monodromy group of the family $\mathcal{Z} \times_B \mathcal{X} \to B$. 
  By Andr\'e's Theorem \ref{AndreInvariantCycle} we get a surjection $\mathfrak{m}(Y) \to \mathfrak{h}(Y)$. Since
  $\mathfrak{m}(Y)$ is an object of $\langle \mathfrak{h}^2(Y) \rangle$
  and since this category is closed under quotients in $\MM$, the proof is reduced to the construction of the $\eta_i$.

  \vspace{1ex}
  Let $\eta_{i,0}: \mathfrak{h}^2(X) (-i+1) \to \mathfrak{h}^{2i}(X)$ 
  be the morphism of motives corres\-ponding to the algebraic homomorphism $\eta_{i,0}$
  in item iii) in the summary at the end of the last section.
  Next, we define $\eta_{i,1}: \uu(-i) \to \mathfrak{h}^{2i}(X)$ 
  as the motivated cohomlogy class $\eta_{i,1} = \varphi_i (v) \in H^{2i}(X,\QQ)$. Finally
  we define
  \begin{equation*}
   \eta_i := \eta_{i,0} \oplus \eta_{i,1} : \mathfrak{g}_i(X) \to \mathfrak{h}^{2i}(X).
  \end{equation*}
  
  Property 1) has been checked in items iii) and iv) at the end of the preceding section.
    
  Property 2) is a direct consequence of item v). There we have seen
  that the Betti realization of $\eta$ is surjective. But the Betti realization is a conservative
  functor. Thus, $\eta$ is surjective in $\mathcal{M}$. This proves a)

\vspace{2ex}
\noindent b) is a direct consequence of \cite[Thm.\ 1.5.1]{An2}. This theorem says that the motive $\mathfrak{h}^2(Y)$ of an irreducible
 symplectic variety $Y$ is an object of $\mathcal{M}(\mathrm{Ab})$, the smallest Tannakian subcategory of the category
 of Andr\'e motives which contains the motives of Abelian varieties.
 The proof of this theorem relies 
 on the Kuga--Satake correspondence. Andr\'e shows that the Kuga--Satake homomorphism $P^2(Y) \hookrightarrow H^2(A \times A, \QQ)$
 is motivated where $P^2(Y)$ is the primitive part of $H^2(Y,\QQ)$ with respect to some polarization and $A$ is a 
 Kuga--Satake variety for $P^2(Y)$. Thus, $\mathfrak{p}^2(Y)$, the motive corresponding to $P^2(Y)$, and hence also
 $\mathfrak{h}^2(Y)$ are objects of $\mathcal{M}(\mathrm{Ab})$.

 \vspace{2ex}
\noindent c) follows from b) and from \cite[Thm.\ 0.6.2]{An1}, which says that all Hodge classes on Abelian varieties
 are motivated. The proof of this theorem uses the deformation principle to reduce first to Abelian varieties with CM, then
 to Weil classes and finally to products of elliptic curves. \qed
\end{section}


\begin{thebibliography}{00000}

\bibitem[An1]{An1} Y.\ Andr\'e, \emph{Pour une th\'eorie inconditionnelle des motifs}, Publ.\ Math.\ I.H.E.S. 
 {\bf 83} (1996), 5-49.
\bibitem[An2]{An2} Y.\ Andr\'e, \emph{On the Shafarevich and Tate conjectures for hyperk\"ahler varieties},
 Math.\ Ann.\ {\bf 305} (1996), 205-248.
\bibitem[Ar]{Ar} D.\ Arapura, \emph{Motivation for Hodge cycles}, Adv.\ Math.\ {\bf 207}, no.\ 2 (2006), 762-781.
\bibitem[B]{Be} A.\ Beauville, \emph{Vari\'et\'es K\"ahl\'eriennes dont la premi\`ere
 classe de Chern est nulle}, J.\ Diff.\ Geom.\ {\bf 18} (1983), 755-782.
\bibitem[J]{Ja} U.\ Jannsen, \emph{Motives, numerical equivalence, and semi-simplicity},
Invent.\ Math.\ {\bf 107}, no.\ 3 (1992), 447-452.
\bibitem[K]{Kl} S.\ Kleiman, \emph{Algebraic cycles and the Weil conjectures}, in: Dix expos\'es sur la cohomologie des sch\'emas, North-Holland, Amsterdam (1968), 359-386.
\bibitem[M1]{Mar1} E.\ Markman, \emph{Generators of the cohomology ring of moduli spaces of 
sheaves on symplectic surfaces}, J.\ reine angew.\ Math.\ {\bf 544} (2002), 61-82.
\bibitem[M2]{Mar2} E.\ Markman, \emph{On the monodromy of moduli spaces of sheaves on K3 surfaces}, J.\ Alg.\ Geom.\
 {\bf 17} (2008), 29-99.
\bibitem[O'G1]{OG} K.\ O'Grady, \emph{The weight-two Hodge structure of moduli spaces of 
sheaves on a K3 surface}, J.\ Alg.\ Geom.\ {\bf 6} (1997), 599-644.
\bibitem[O'G2]{OG1} K.\ O'Grady, \emph{Desingularized moduli spaces of sheaves on a K3}, 
J.\ Reine Angew.\ Math.\ {\bf 512} (1999), 49-117.
\bibitem[O'G3]{OG2} K.\ O'Grady, \emph{A new six-dimensional irreducible symplectic variety},
J.\ Alg.\ Geom.\ {\bf 12}, no.\ 3 (2003), 435-505.
\end{thebibliography}
\end{document}